Collaborative Preferences for Learning Mathematics: A Scale Validation Study

Sang Hyun Kim
University of Auckland

Tanya Evans
University of Auckland

*Collaboration within mathematics has been established as being effective in providing students with crucial opportunities to develop critical thinking, effective communication, and teamwork skills. By engaging in group problem-solving and shared learning experiences, students may gain deeper insights into mathematical concepts and learn to approach challenges from multiple perspectives. However, there remains a need for a reliable instrument to capture students' preferences for collaboration. This study aims to develop and validate the Collaborative Preferences for Learning Mathematics (CPLM) scale to measure student preferences for collaboration. Exploratory factor analysis revealed a single-factor structure, and a confirmatory factor analysis conducted with a separate sample (N = 243) demonstrated a good model fit. Further testing established the scale's strong invariance, confirming its ability to reliably measure collaborative preferences over time. The CPLM offers a valid and reliable way to capture student preferences regarding collaborative learning in mathematics.*

Undergraduate students are offered various opportunities to engage with mathematics, including a diverse range of pedagogical approaches, assessments, and resources. It is important to acknowledge that the backgrounds, experiences, and values students bring influence their preferences for how they engage with mathematics, particularly with their peers. Understanding and accounting for these preferences enables educators to tailor their teaching strategies, create more effective group dynamics, and foster inclusive environments that accommodate diverse learner needs. However, the extent to which students prefer collaboration or individual work, especially in the context of mathematics, remains underexplored, and there is a notable lack of systematic, quantitative tools to capture these preferences on a large scale. To address this gap, this study sought to develop and validate the Collaborative Preferences for Learning Mathematics (CPLM) scale, providing a robust measure to better understand these preferences. In the following sections, We review the literature on student perceptions of collaborative practices in mathematics, outline the study's context, present the results, and discuss the implications.

**Literature Review**

Collaboration is widely established in the literature as an efficacious approach to learning. Grounded in social learning theories such as Bandura's social-cognitive theory (1997), collaborative settings allow students to learn vicariously, increasing self-efficacy. Collaborative learning positively affects achievement and attitudes towards learning across various disciplines (Springer et al., 1999; Tullis & Goldstone, 2020; Xu et al., 2001). Specifically in mathematics, it has been argued that collaboration positively shapes students' learning experiences and could be used as an approach to foster motivation, confidence, and achievement for some students (Cuneo, 2008; Duah et al., 2014). However, achievement is not directly predicted by the time spent studying or the proportion spent studying collaboratively (Alcock et al., 2020).

Despite these general benefits, not all students have the same inclination to work with others (Johnson & Johnson, 1995); collaborative and individual learning preferences vary significantly. For instance, gender and age have been shown to impact this, with girls showing stronger preferences for cooperative learning than boys, who tend to favour competitive and individualised learning approaches. Additionally, preferences for cooperative and competitive learning have been reported to increase with school age (Owens, 1985). Personality traits may also be a predicting factor in students' willingness to learn collaboratively. One study reported that extroversion was significantly negatively correlated with a preference for independent study, and neuroticism, a measure of emotional stability, was negatively correlated with preferences for social modes of learning (Chamorro-Premuzic et al., 2007). Furnham et al. (2008) reported that group work was preferred by younger extroverted students with low conscientiousness and low general knowledge. Furthermore, these preferences and attitudes may also vary across disciplines (Gottschall & García-Bayonas, 2008), where disciplinary norms play a significant role. In mathematics, an emphasis on individual problem-solving, commonly used teacher-centered instruction (Dietrich & Evans, 2022; Melhuish et al., 2022), and assessments that focus on procedural knowledge may limit opportunities for collaboration and reinforce a perception that mathematics is best learned independently rather than collaboratively.

**Student Perspectives of Collaborative Learning**

Beyond individual factors, there is considerable nuance in students' experiences that are shaped by the social context in which they operate. In broader educational contexts, students' collaboration preferences often reflect their priorities and inclinations. Some students may choose to study alone to minimize distractions, while others may prefer group study to expose themselves to diverse perspectives. Many of the key difficulties of group work identified in the broader education literature tend to highlight the importance of equal contribution from group members (Feichtner & Davis, 1984; Salomon & Globerson, 1989; Tucker & Abbasi, 2016). This issue can be especially problematic when group work is assessed and has also been reported in the mathematics education literature.

Within the context of undergraduate mathematics, one study that interviewed three cohorts of students across different universities revealed that while students generally had positive attitudes towards group work, many still expressed a preference for learning individually in mathematics (MacBean et al., 2004). Within their study, group work, which was largely experienced in informal settings outside of class, provided a sense of camaraderie between students. These casual meetings offered the time for students to discuss their collective difficulties. Additionally, the student-driven nature of these informal interactions allows students to have the freedom to use the time as they desire with peers with whom they get along. Group work done in a formal capacity, such as during lectures or for assessments, can provide similar benefits but also expose numerous drawbacks. One such drawback is the occurrence of social loafing, which refers to "the reduction in motivation and effort when individuals work collectively compared with when they work individually or coactively" (Karau & Williams, 1993, p. 681), is more salient in larger groups than in smaller groups (North et al., 2000). Working with unfamiliar group members was another issue mentioned by students who must work to establish new relationships whilst navigating mathematical challenges (MacBean et al., 2004).

D'Souza and Wood (2003) corroborated the findings of MacBean et al. (2004): students perceived group work as a valuable experience with potential cognitive benefits and saw it as an opportunity to meet new people and increase self-esteem. However, they also identified significant challenges, including the presence of less outspoken and unprepared students who were reluctant to participate, as well as individuals who did not contribute equally. The latter aspect might be advantageous for disengaged or weaker students who benefit from others' efforts. This observation prompts a nuanced discussion about the balance between supporting weaker group members while managing unequal workloads and maintaining fairness within collaborative environments.

Further research exploring student attitudes toward group work revealed generally positive views, with most students believing that collaboration enhances the effectiveness of mathematics teaching (Sofroniou & Poutos, 2016). However, like the previous studies, Sofroniou and Poutos (2016) reported concerns about unproductive group dynamics that were salient among students. Additionally, some students even recognized that group work "can slow down the lesson" (p. 7).

While collaborating in mathematics can be difficult for students, little is known why this might be the case. One reason could be that the strong tendencies of students to compare themselves with peers may hinder their ability to collaborate effectively in mathematics (Campbell & Yeo, 2023). When reflecting on their actions during collaborative activities, students demonstrate an acute awareness of various elements of their interactions – such as how often they speak and the overall group rapport – and the underlying factors that influence these dynamics, including confidence and personality (Campbell & Yeo, 2023).

## Research Aims

The aim of this study was two-fold: (1) to develop and validate the CPLM scale, an instrument intended to measure general student preferences towards learning mathematics, and (2) to establish measurement invariance of the CPLM.

**Rationale**

Given the challenges and potential benefits of collaborative learning in mathematics, it is essential to better understand and gauge students' preferences for learning in this context. While collaborative learning has been explored in other subjects, there remains a need for research specifically focused on mathematics education. Much empirical and theoretical research highlights the various ways that working with others in educational settings can impact student learning; researchers should maintain a keen interest in students' preferences in these educational settings.

Firstly, this exploration of student voice indicates what students, a significant stakeholder in education, are experiencing and valuing in their learning. Kouros et al. (2006) describe students' attitudes towards learning as a reflection of the quality of their experiences. Understanding their preferences for individual versus collaborative learning in mathematics can provide valuable insights into how instructional approaches can be tailored to meet their needs better and support these experiences. Secondly, in a rapidly evolving world where social interaction preferences are continually shifting, having accurate methods for evaluating students' desires and needs can be invaluable for our understanding of the nuanced views from the students' perspectives, which has the potential to enhance overall student engagement and satisfaction. Lastly, even though student views represent just one overlooked dimension of mathematical learning, it can help identify whether the practices and policies that are implemented are aligned with their wants.

This research does not aim to advocate for or against any specific instructional method or any particular attitude towards students' collaborative preferences. Instead, it offers a way to get a snapshot of undergraduate students' preferences in collaborative settings within a mathematics context. While much existing research has focused on samples from different educational levels, our study seeks to explore straightforward and accessible methods to gauge student preferences in learning. Specifically, this work captures undergraduate student experiences in a more generalized and context-independent manner, providing a broad overview of current student preferences. Developing a validated scale to measure student preferences for collaborative learning in mathematics can provide educators with insights into how instructional strategies might be adjusted to better align with student needs and promote engagement in this traditionally individualistic discipline.

## Methodology

**Collaborative Preferences for Learning Mathematics (CPLM)**
   This short scale contains five items that prompt students to consider their preferences for collaboration in mathematics. The items gauge different contexts and settings where these preferences may vary, including when learning for high- and low-stakes assessments and in different stages of their mathematical competency (being exposed to new concepts). All items were responded to on a slider scale from 0 to 100. The prompt and items in the scale are shown below.
   Prompt: "Consider yourself learning mathematics. By moving the slider, state the extent to which you prefer to do it from 0 = Individually to 100 = Collaboratively for each of the items:"
  1. What is the most effective way for you to learn mathematics? (CPLM_1)
  2. What is the best way for you to make sense of mathematics? (CPLM_2)
  3. What is the most effective way for you to study for high-stakes maths assessments (e.g., exams)? (CPLM_3)
  4. In what social setting do you prefer to be exposed to novel concepts? (CPLM_4)
  5. In what social setting do you prefer to engage in problem-solving in low-stakes assessment (e.g., homework, practice exercises)? (CPLM_5)

**Analysis**
   **Exploratory Factor Analysis (EFA).** To validate the CPLM, we conducted an EFA to confirm the underlying structure of the five items (observed variables) and the latent variable (CPLM). Exploratory factor analysis was performed using principal axis factoring as the extraction method, which is more appropriate for non-normal data (Costello & Osborne, 2005). We assessed the appropriateness of the data for factor analysis using the Kaiser-Meyer-Olkin (KMO) measure of sampling adequacy and Bartlett's test of sphericity. A KMO value above the .50 threshold and the rejection of the null hypothesis for Bartlett's test were considered indicators that the data was appropriate for this analysis (Dziuban & Shirkey, 1974). Cronbach's alpha values greater than .70 were considered satisfactory (Taber, 2018). The EFA was conducted on SPSS (Version 29).
   **Confirmatory Factor Analysis (CFA).** We conducted a CFA to validate and confirm the factor structure identified in the EFA. A 'good' model fit was defined by $\chi^2/df < 3.0$, TLI > .95, CFI > .95, RMSEA < .06, and SRMR < .05. An 'acceptable' fit was indicated by $\chi^2/df < 5.0$, TLI > .90, CFI > .90, RMSEA < .08, and SRMR < .08 (Hu & Bentler, 1999). We assessed convergent validity by evaluating whether the Average Variance Extracted (AVE) exceeded .50 and Composite Reliability (CR) was greater than .70. AMOS (Version 29) was used for the CFA and invariance testing.
   **Multi-Group Confirmatory Factor Analysis (MGCFA).** We conducted an MGCFA to test the measurement invariance of the CPLM over three different time points. We compared changes in model fit indices to evaluate the invariance, focusing on the Comparative Fit Index (CFI). A change in CFI of less than 0.01 between models was considered indicative of good measurement invariance (Cheung & Rensvold, 2002).

**Data**

In this study, we utilized self-report data collected as part of a larger project. Participants rated their preferences on a 0 to 100 scale using five slider-based, close-ended items. The questionnaires were distributed via Qualtrics, an online survey platform. Two data samples were collected. The first sample consisted of undergraduate students recruited through Prolific, an online participant recruitment platform. These participants were current undergraduates from the UK or USA, studying STEM or commerce-related majors. After data cleaning, 273 valid responses were used for EFA. The second dataset was collected during the second semester of 2023 from a second-year undergraduate service mathematics course at a large university in New Zealand. This sample involved repeated measures from 224 students across three time points. Data from the first time point (N = 243) was used for CFA, and the repeated measures across the three time points were used for the multi-group confirmatory factor analysis (MGCFA) to test for measurement invariance over time.

## Results

**Exploratory Factor Analysis.** The inter-item correlations of the CPLM ranged from .209 to .619. Bartlett's Test of Sphericity was significant, $\chi^2(10) = 322.69$, $p < .001$, and the KMO measure of sampling adequacy was .759, both indicating that the data were suitable for factor analysis. The results of the EFA yielded a single factor with an eigenvalue of 2.005, accounting for 40.09% of the variance. All five items loaded adequately to strongly to the factor (.476 to .863) (Costello & Osborne, 2005) and had communalities ranging from .227 to .744 (see Table 1). The five items in the factor show an acceptable degree of internal consistency, with a Cronbach's alpha value of .744.

*Table 1. Factor loadings and communalities for each item from the EFA.*

| CPLM | | |
|---|---|---|
| Item | Factor Loading | Communalities |
| CPLM_1 | .863 | .744 |
| CPLM_2 | .734 | .539 |
| CPLM_3 | .507 | .257 |
| CPLM_4 | .488 | .238 |
| CPLM_5 | .476 | .227 |

**Confirmatory Factor Analysis.** The results of the CFA showed that the one-factor model demonstrated a good fit to the data: $\chi^2/df = 1.601$, TLI = .988, CFI = .994, RMSEA = .050, and SRMR = .0230 (see Figure 1). The factor loadings ranged from .59 to .87. The Cronbach's alpha for the factor was .848, indicating good internal consistency among the items. Additionally, the Average Variance Extracted (AVE) was 0.54 and the Composite Reliability (CR) was 0.85, indicating that convergent validity was achieved. This evidence provided a strong basis for proceeding with the MGCFA to test for measurement invariance.

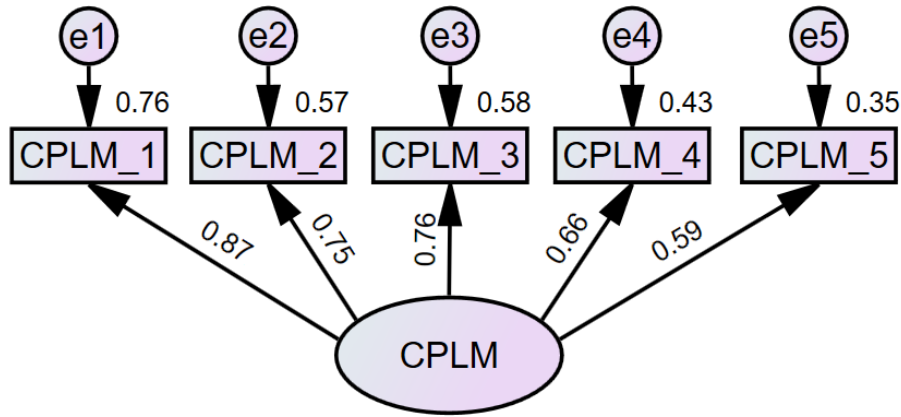

*Figure 1. One-factor model for the CPLM.*

**Invariance Testing.** The changes in CFI across the models remained within 0.01 of each model, except for the residual model, which had a change in CFI of .043 (see Table 2). From this, we can conclude that the model demonstrates strong invariance (configural, metric, scalar) but lacks strict invariance (configural, metric, scalar, and residual). Establishing strong invariance is generally regarded as adequate for determining model equivalence, allowing for valid comparisons of latent factor means across time (Meredith, 1993). Strict invariance, while rarely achieved in practice due to its stringent requirement of identical reliability across groups, is typically not considered necessary unless the research question specifically involves comparing observed variances or covariances across groups.

*Table 2. Measurement invariance of the CPLM over a semester (N = 224).*

| Model | CFI | ΔCFI | $\chi^2/df$ | RMSEA |
| --- | --- | --- | --- | --- |
| Configural (Unconstrained) | .993 | - | 1.804 | 0.035 |
| Metric | .990 | .003 | 1.803 | 0.035 |
| Scalar | .991 | .001 | 1.671 | 0.032 |
| Residual | .948 | .043 | 3.742 | 0.064 |

### Discussion

This study aimed to develop and validate an instrument to measure students' general preferences for collaboration within mathematics. The results suggested that the CPLM scale is reliable and valid for measuring these preferences. EFA identified a single-factor structure, while CFA indicated the model demonstrated a good fit. Furthermore, MGCFA provided evidence of strong invariance, confirming that the scale reliably measures the same construct across different time points.

One limitation of the CPLM is that, while it was validated in an undergraduate context, its application in other settings will necessitate further investigation and validation. Future research could expand the scale to include a broader range of collaborative experiences across different contexts, thereby enhancing its comprehensiveness and applicability.

Additionally, incorporating qualitative data, such as interviews or observations, could provide a richer and more nuanced understanding of students' collaborative preferences and the underlying factors that shape them. These insights could inform the refinement of existing tools and frameworks, ensuring they better capture the complexities of collaboration in educational settings. Such approaches would also create opportunities to explore how these preferences impact group dynamics, engagement, and overall learning experiences.

Despite the limitations, we believe this instrument offers a practical and efficient way to assess collaborative learning preferences on a large scale. Its ability to quickly gauge these preferences provides educators with valuable insight into their students' preferences towards collaboration. Such insight would enable instructors to anticipate and plan for diverse preferences within their classrooms, allowing them to better structure group activities and allocate resources. The CPLM scale is not intended to advocate for pedagogical practices that exclusively align with student preferences. Therefore, educators should use CPLM results as one of many tools to inform their teaching strategies rather than as prescriptive guidelines.

From a research perspective, the CPLM scale offers a quantitative tool to assess these preferences systematically, allowing researchers to incorporate student learning preferences into studies, which could serve as a useful variable to account for alongside other factors in future research. This integration can deepen our understanding of how collaborative dynamics influence learning outcomes and student engagement. Longitudinal studies, in particular, could shed light on how these preferences evolve over time or in response to targeted instructional interventions. In summary, the CPLM scale contributes to advancing educational research by offering a reliable tool for understanding collaborative learning preferences in mathematics.

## References

Alcock, L., Hernandez-Martinez, P., Godwin Patel, A., & Sirl, D. (2020). Study Habits and Attainment in Undergraduate Mathematics: A Social Network Analysis. *Journal for research in mathematics education*, *51*(1), 26-49. https://doi.org/10.5951/jresematheduc.2019.0006

Bandura, A. (1997). *Self-efficacy : the exercise of control*. New York : W.H. Freeman c1997.

Campbell, T. G., & Yeo, S. (2023). Student noticing of collaborative practices: exploring how college students notice during small group interactions in math. *Educational Studies in Mathematics*, *113*(3), 405-423. https://doi.org/10.1007/s10649-023-10206-3

Chamorro-Premuzic, T., Furnham, A., & Lewis, M. (2007). Personality and approaches to learning predict preference for different teaching methods. *Learning and Individual Differences*, *17*(3), 241-250. https://doi.org/https://doi.org/10.1016/j.lindif.2006.12.001

Cheung, G. W., & Rensvold, R. B. (2002). Evaluating Goodness-of-Fit Indexes for Testing Measurement Invariance. *Structural Equation Modeling: A Multidisciplinary Journal*, *9*(2), 233-255. https://doi.org/10.1207/S15328007SEM0902_5

Costello, A. B., & Osborne, J. W. (2005). Best practices in exploratory factor analysis: Four recommendations for getting the most from your analysis. *Practical assessment, research & evaluation*, *10*(7), 7.

Cuneo, A. (2008). *Examining the effects of collaborative learning on performance in undergraduate mathematics* (Publication Number 3290667) [Ph.D., Capella University]. ProQuest Dissertations & Theses Global. United States -- Minnesota.

D'Souza, S. M., & Wood, L. N. (2003). Tertiary students' views about group work in mathematics. Proceedings of the Educational Research, Risks and Dilemmas—New Zealand Association for Research in Education (NZARE) and Australian Association for Research in Education (AARE) Joint Conference, The University of Auckland, Auckland, New Zealand,

Dietrich, H., & Evans, T. (2022). Traditional lectures versus active learning – A false dichotomy? *STEM Education*, *2*(4), 275-292. https://doi.org/10.3934/steme.2022017

Duah, F., Croft, T., & Inglis, M. (2014). Can peer assisted learning be effective in undergraduate mathematics? *International Journal of Mathematical Education in Science and Technology*, *45*(4), 552-565. https://doi.org/10.1080/0020739X.2013.855329

Dziuban, C. D., & Shirkey, E. C. (1974). When is a correlation matrix appropriate for factor analysis? Some decision rules. *Psychological Bulletin*, *81*(6), 358-361. https://doi.org/10.1037/h0036316

Feichtner, S. B., & Davis, E. A. (1984). Why Some Groups Fail: a Survey of Students' Experiences with Learning Groups. *Organizational Behavior Teaching Review*, *9*(4), 58-73. https://doi.org/10.1177/105256298400900409

Furnham, A., Christopher, A., Garwood, J., & Martin, N. G. (2008). Ability, demography, learning style, and personality trait correlates of student preference for assessment method. *Educational Psychology*, *28*(1), 15-27. https://doi.org/10.1080/01443410701369138

Gottschall, H., & García-Bayonas, M. (2008). Student Attitudes Towards Group Work Among Undergraduates in Business Administration, Education and Mathematics. *Educational Research Quarterly*, *32*(1), 3-28.

Hu, L. t., & Bentler, P. M. (1999). Cutoff criteria for fit indexes in covariance structure analysis: Conventional criteria versus new alternatives. *Structural Equation Modeling: A Multidisciplinary Journal*, *6*(1), 1-55. https://doi.org/10.1080/10705519909540118

Johnson, D. W., & Johnson, R. T. (1995). Cooperative Learning and Nonacademic Outcomes of Schooling. In J. E. Pedersen & A. D. Digby (Eds.), *Secondary Schools and Cooperative Learning* (pp. 81-150). Garland.

Karau, S. J., & Williams, K. D. (1993). Social Loafing: A Meta-Analytic Review and Theoretical Integration. *Journal of personality and social psychology*, *65*(4), 681-706. https://doi.org/10.1037/0022-3514.65.4.681

Kouros, C., Abrami, P. C., Glashan, A., & Wade, A. (2006). How do students really feel about working in small groups? The role of student attitudes and behaviors in cooperative classroom settings. annual meeting of the American Educational Research Association, San Francisco, California,

MacBean, J., Graham, T., & Sangwin, C. (2004). Group work in mathematics: a survey of students' experiences and attitudes. *Teaching Mathematics and its Applications: An International Journal of the IMA*, *23*(2), 49-68. https://doi.org/10.1093/teamat/23.2.49

Melhuish, K., Fukawa-Connelly, T., Dawkins, P. C., Woods, C., & Weber, K. (2022). Collegiate mathematics teaching in proof-based courses: What we now know and what we have yet to learn. *The Journal of Mathematical Behavior*, *67*, 100986. https://doi.org/https://doi.org/10.1016/j.jmathb.2022.100986

Meredith, W. (1993). Measurement invariance, factor analysis and factorial invariance. *Psychometrika*, *58*(4), 525-543. https://doi.org/10.1007/BF02294825